\documentclass{ifacconf}

\usepackage{graphicx}      
\usepackage{natbib}        

\usepackage[utf8]{inputenc}
\usepackage{graphics} 
\usepackage{epsfig} 
\usepackage{epstopdf}
\usepackage{amsmath} 
\usepackage{amssymb}  
\usepackage{amsfonts}

\newcommand{\dt}{\mathrm{dt}}

\usepackage{xspace}
\usepackage{soul}
\usepackage{mathtools}
\usepackage{makecell}
\usepackage{cite}
\usepackage{nicefrac}
\usepackage{float}
\usepackage{wrapfig}
\usepackage{xcolor}

\newcommand{\ra}{\rightarrow}

\newcommand{\Ra}{\Rightarrow}

\newtheorem{dfn}{Definition}
\newtheorem{asm}{Assumption}

\newcommand{\ie}{\unskip, i.\,e.,\xspace}
\newcommand{\eg}{\unskip, e.\,g.,\xspace}

\newcommand{\sut}{\text{s.\,t.\,}}

\newcommand{\N}{\ensuremath{\mathbb{N}}}

\newcommand{\R}{\ensuremath{\mathbb{R}}}

\newcommand{\X}{\ensuremath{\mathbb{X}}}

\newcommand{\U}{\ensuremath{\mathbb{U}}}

\newcommand{\B}{\ensuremath{\mathcal{B}}}

\newcommand{\eps}{\ensuremath{\varepsilon}}

\newcommand{\ball}{\ensuremath{\mathcal B}}
\newcommand{\D}{\ensuremath{\mathcal{D}}}		
\newcommand{\K}{\ensuremath{\mathcal{K}}\xspace}		

\newcommand{\nrm}[1]{\left\lVert#1\right\rVert}

\newcommand{\scal}[1]{\left\langle#1\right\rangle}

\definecolor{dgreen}{rgb}{0.0, 0.5, 0.0}

\newcommand{\spc}{\ensuremath{\,\,}}	

\makeatletter
\newcommand{\subalign}[1]{%
	\vcenter{%
		\Let@ \restore@math@cr \default@tag
		\baselineskip\fontdimen10 \scriptfont\tw@
		\advance\baselineskip\fontdimen12 \scriptfont\tw@
		\lineskip\thr@@\fontdimen8 \scriptfont\thr@@
		\lineskiplimit\lineskip
		\ialign{\hfil$\m@th\scriptstyle##$&$\m@th\scriptstyle{}##$\crcr
			#1\crcr
		}%
	}
}

\newcommand{\norm}[1]{\left\lVert#1\right\rVert}

\begin{document}
	\setlength{\belowdisplayskip}{2.5pt}
\begin{frontmatter}

\title{A reinforcement learning method with closed-loop stability guarantee} 


\author[First]{Pavel Osinenko} 
\author[First]{Lukas Beckenbach} 
\author[First]{Thomas G\"{o}hrt}
\author[First]{Stefan Streif}

\address[First]{Technische Universit\"at Chemnitz, Automatic Control and System Dynamics Laboratory, Germany (e-mail: p.osinenko@gmail.com, \{lukas.beckenbach,thomas.goehrt,stefan.streif\}@etit.tu-chemnitz.de)}


\begin{abstract}                
Reinforcement learning (RL) in the context of control systems offers wide possibilities of controller adaptation.
Given an infinite-horizon cost function, the so-called critic of RL approximates it with a neural net and sends this information to the controller (called ``actor'').
However, the issue of closed-loop stability under an RL-method is still not fully addressed.
Since the critic delivers merely an approximation to the value function of the corresponding infinite-horizon problem, no guarantee can be given in general as to whether the actor's actions stabilize the system.
Different approaches to this issue exist.
The current work offers a particular one, which, starting with a (not necessarily smooth) control Lyapunov function (CLF), derives an online RL-scheme in such a way that practical semi-global stability property of the closed-loop can be established.
The approach logically continues the work of the authors on parameterized controllers and Lyapunov-like constraints for RL, whereas the CLF now appears merely in one of the constraints of the control scheme.
The analysis of the closed-loop behavior is done in a sample-and-hold (SH) manner thus offering a certain insight into the digital realization.
The case study with a non-holonomic integrator shows the capabilities of the derived method to optimize the given cost function compared to a nominal stabilizing controller.
\end{abstract}
\begin{keyword}
Reinforcement learning control, Stability of nonlinear systems, Lyapunov methods 
\end{keyword}

\end{frontmatter}

\section{Introduction}\label{sec:intro}

Consider a general nonlinear dynamical system
\begin{equation}	\label{eq:sys}
\dot x = f(x,u),
\end{equation}
where $x \in \R^n$ is the state, $u \in \R^m$ is the input, $f : \R^n \times \R^m \to \R^n$ is the 
dynamics model.
Further, consider the following infinite-horizon (IH) cost function:
\begin{equation}	\label{eq:IH-cost}
J[\kappa](x_0) := \int \limits_0^{\infty} \rho(x(t), \kappa(x(t))) \; \dt, \quad x(0) = x_0,
\end{equation}
where $\rho : \R^n \times \R^m \to \R_{\geq0}$ denotes the \textit{reward} function and $\kappa : \R^n \to \R^m$ is a control policy.
The function $J^*(x_0) := \min_\kappa J[\kappa](x_0), \forall x_0$ is called the \textit{value function}.
By the Bellman's optimality principle, it satisfies the Hamilton-Bellman-Jacobi equation:
\begin{equation}	\label{eq:HJB}
\dot J^*(x) + \min_u \{\nabla J^* f(x,u) + \rho(x,u)\} = 0, \forall x.
\end{equation}
Dynamic programming (DP) takes \eqref{eq:HJB} as the basis, discretized (in a compact domain of) the state space, and computes and approximation to $J^*$ in an iterative manner \citep{Liu2014-PI-ADP-conv-prf,Wei2016}.
The \textit{curse of dimensionality} is what prevents application of DP in the dynamical context.
One particular way to overcome this issue is to use parameterized function approximators $\hat J (x,\theta)= \scal{\theta, \varphi(x)}$ with a finite number of parameters, e.g. neural nets. 
Here, $\theta$ is the hidden layer weight vector and $\varphi$ is the activation function of the net.
Thus, the iterations are performed over the parameters $\theta$.
Roughly, the idea reads:
\begin{equation}	\label{eq:RL-idea}
\begin{aligned}
\text{Step 1:} \spc & \theta_{\text{new}} := \min_{\theta} \; \{\rho + \Delta \hat J\} & (\text{Critic}),\\
\text{Step 2:} \spc & u_{\text{new}} := \min_u \;  \{\rho + \hat J\} & (\text{Actor}),\\
\end{aligned}
\end{equation}
where $\rho + \Delta \hat J$ represents the \textit{Bellman error} \ie a metric, which describes the goodness of $\hat J$ as an approximation to $J^*$ based on the HJB.
There is a variety of RL methods, 
thus it is virtually impossible to comprehensively overview them \citep[the reader may refer \eg to][]{Bertsekas2017,Sutton2018,Recht2019}. 
However, it is worthwhile to categorize some in terms of how they tackle the issue of closed-loop stability, since direct application of \eqref{eq:RL-idea} does not necessarily 
give any guarantees.

There are methods that: (a) are heavily based on DP principles \citep{Heydari2014a,Wei2016}, 
(b) concentrate solely on neural net weight learning \citep{Sokolov2015,Zhang2011,Mu2017}, 
(c) start with sufficiently good initial data \citep{Jiang2015,Gao2017}, 
(d) restrict to linear systems \citep{Bian2016,Gao2016}. 
The first category entails iterations over (a subset) of the state space, (b) require long off-line learning phases and do not take into account closed-loop stability, (c) puts the burden of fine initialization onto the user.
In general, there is oftentimes a dilemma: RL pursues optimality of DP, but often lacks stability guarantees, whereas some nominal stabilizing controller is not concerned about optimality in the sense of minimizing \eqref{eq:IH-cost}. 
The relations between optimal and stabilizing controllers were well described in \citep{Primbs1999}. 
It seems a certain trade-off is required to tackle optimality and closed-loop stability simultaneously.

The \textbf{contribution} of the current work is to offer an RL-method which does address closed-loop stability.
It is based on an initial stabilizability information, specifically, in the form of a (not necessarily smooth) CLF.
The justification of such an assumption is as follows.
Every existing RL approach requires at least stabilizability of the system.
Stabilizability implies in turn existence of a CLF by a converse result.
It is suggested to constrain the RL-method accordingly.
Similar philosophy was pursued in the previous work of the authors \citep{Beckenbach2018MPC,2019_goeht_adp_lyap}. 
However, the current work greatly generalizes the previous derivations.
First, the assumed CLF needs not to be smooth, as it is in the general case \citep{Clarke1997-stabilization}. 
Secondly, state convergence shown in this work is provided in the sense of practical stabilizability instead of just ultimate boundedness used some literature \citep[see \eg][]{Vamvoudakis2010}.
The new algorithm is thus suggested in a sample-and-hold setting (SH) which gives insight into the digital realization.
In particular, the actor and critic are now merged, and the ``actor-critic'' optimization is performed at discrete time samples.
Roughly, the method reads:
\begin{equation}	\label{eq:RL-new}
\begin{aligned}
(u_{\text{new}}, w_{\text{new}}) \spc & := \min_{(u,w)} \; \mathbb J (x,u,w) & (\text{Actor-Critic}),\\
\sut \spc & \Delta_{\text{inter-sample}} \hat J < 0 & (\text{Constraints}),\\
& \Delta_{\text{sample-to-sample}} \hat J < 0
\end{aligned}
\end{equation}
where $\mathbb J$ is a cost function related to the Bellman error.
Notice here the requirement of inter-sample and sample-to-sample decay of the critic $\hat J$.
The actual algorithm, which is presented in Section \ref{sec:main-results}, does not pose the constraints so literally -- there are certain relaxation terms.
The case study with a non-holonomic integrator demonstrates the worthiness of \eqref{eq:RL-new} in Section \ref{sec:case-study}. 

\textit{Notation.} A closed ball of radius $R>0$ centered at the origin is denoted $\B_R$. 
A continuous function $\alpha:[0,a) \ra [0,\infty)$ is said to belong to class $\mathcal{K}$, if it is strictly increasing and $\alpha(0) = 0$. 
It is said to belong to $\mathcal{K}_{\infty}$ if $a=\infty$ and $\alpha(r) \ra \infty$ as $r \ra \infty$.  


\section{Preliminaries}\label{sec:preliminaries}

As mentioned above, the suggested RL-method will be considered in a SH setting.
Such a setting means applying constant controls throughout sampling periods of some time $\delta>0$, in which the system 
is governed by
\begin{align} \label{eq:sys-SH}
\begin{split}
&\dot{x} = f(x, u_k), \quad x(0), \\
&t \in [k \delta,(k+1) \delta), \; u_k = \kappa(x(k \delta)), \; k \in \N_0,
\end{split}
\end{align} 
where $u_k \in \U \subset \R^m$ are input constraints. 
It is assumed that the dynamics model $f$ is locally Lipschitz in $x$ for any $u \in \U$. 
In the following denote $x(k \delta) =: x_k$. 
For any $k \in \N_0$ and $\delta > 0$, the state $x^{u_k}(t)$ at $t \geq k \delta$ under $u_k$ is defined as
\begin{align}
x^{u_k}(t) := x_k + \int_{k \delta}^{t} f(x(\tau), u_k) \; \text{d} \tau.
\end{align}
For $t = (k+1) \delta$, denote $x_{k+1}^{u_k} := x^{u_k}((k+1)\delta)$. 
The corresponding trajectory of \eqref{eq:sys-SH} under the SH-mode input will also be called \textit{SH-trajectory}. 
Consider the following standard

\begin{dfn} \label{def:pract-stab}
	A control policy $\kappa(\cdot)$ is said to practically semi-globally stabilize \eqref{eq:sys} if, given $0<r<R<\infty$, there exists a $\bar \delta > 0$ \sut any SH-trajectory $x(t)$ with a sampling time $\delta \le \bar \delta$, starting in $\mathcal{B}_R$ is bounded, enters $\mathcal{B}_r$ after a time $T$, which depends uniformly on $R,r$, and remains there for all $t \geq T$.
\end{dfn}
In the following, for brevity, the wording ``semi-globally'' is omitted. 
In the light of Definition \ref{def:pract-stab}, the balls $\mathcal{B}_R$ and $\B_r$ are denoted the \emph{starting} and \emph{target ball}, respectively. 
Recall further the following 

\begin{dfn} \label{def:GLDD}
	For a locally Lipschitz function $V:\R^n \ra \R$ and a $v \in \R^n$, the \emph{generalized lower directional derivative} (GLDD) of $V$ in the direction of $v$ at $x$ is defined as \citep{Sontag1995}
	\begin{align} \label{eq:GLDD}
	\D_{v} V(x) := \liminf_{\tau \ra 0_{+}} \frac{1}{\tau} \left( V(x + \tau v) - V(x) \right).
	\end{align}
\end{dfn}

The following is a stabilizability assumption:

\begin{asm} \label{asm:CLF}
	There exists a locally Lipschitz continuous, positive-definite function $V:\R^n \ra \R_{\geq 0}$, a continuous positive definite function $w:\R^n \ra \R_{\geq 0}$ and $\alpha_{1,2} \in \mathcal{K}_{\infty}$ \sut for any compact $\X \subset \R^n$, there exists a compact set $\U_\X \subseteq \U$ and it holds that, for any $x \in \X$, 
	\begin{enumerate}
		\item[i)] $V$ has a decay rate satisfying
		\begin{align} \label{eq:CLF-decay}
		\inf_{u \in \U_\X} \; \mathcal{D}_{f(x,u)} V(x) \leq - w (x),
		\end{align}
		\item[ii)] $\alpha_1(\|x\|) \leq V(x) \leq \alpha_2(\|x\|)$,
	\end{enumerate}
\end{asm} 
The pair $(V,w)$ is also referred to as a \textit{CLF-pair}.
\begin{rem} \label{rem:pract-stab-policy}
	Under the existence of a CLF as per Assumption \ref{asm:CLF}, practical stabilization in the sense of SH can be realized as follows:
	given a CLF-pair $(V,w)$ and balls with radii $0<r<R<\infty$, there exists a $\bar{\delta}>0$ \sut for any $0 < \delta \le \bar \delta$,
	there is a (possibly discontinuous) map $\mu:\R^n \ra \U$ \sut the SH-trajectory of \eqref{eq:sys-SH} with the sampling period $\delta$ satisfies:
	\begin{align} \label{eq:decay-V-sample}
	V\left(x_{k+1}^{\mu(x_k)}\right) - V(x_k) \leq - \frac{\delta }{2}w(x_k).
	\end{align}
	There are various methods of deriving a SH realization of $\mu$ (refer \eg to \citet{Clarke1997-stabilization,Braun2017}).
	The meaning of the last displayed inter-sample decay condition is that, using \eqref{eq:CLF-decay}, one may calculate such control actions $\mu(x_k)$ at the sampling nodes, that at least half the decay is retained (the relaxation comes from the inter-sample behavior).
\end{rem}
Now, address the actor-critic setup of the paper.
First, the critic 
\begin{align} \label{eq:Jhat}
\hat{J}:\R^n \times \R^p \ra \R_{\geq 0}:(x,\theta) \mapsto \langle \theta,\varphi(x) \rangle,
\end{align}
can be regarded as a neural net, consisting of the hidden layer weights $\theta \in \R^p$ and a locally Lipschitz continuous activation function $\varphi:\R^n \ra \R^p$. 
On one hand, given $\hat{J}$, it holds that on any compact $\X \subset \R^n$, $\hat{J}(x,\theta) \leq \|\theta\| L_{\varphi} \|x\|$, where $L_{\varphi}$ is the corresponding local Lipschitz constant on $\X$.
On the other hand, let the activation function satisfy the following condition: there exists $\underline l \in \K$ \sut for any $x \in \R^n$ and $\theta \in \R^p$, it holds that $\langle \theta,\varphi(x) \rangle \ge \underline l (\nrm{x}) \cdot \nrm{\theta}$.
%



Let the following assumption on the activation function of \eqref{eq:Jhat} and the CLF $V$ as per Assumption \ref{asm:CLF} hold:
\begin{asm} \label{asm:struct-equiv}
	There exists $\theta^{\#} \in \R^p$ \sut $V(x) = \hat{J}(x,\theta^{\#})$, for all $x \in \R^n$. 
\end{asm} 
Assumption \ref{asm:struct-equiv} states that the structure of $\varphi$, which is a designer's choice, be ``rich'' enough to structurally capture $V$, which is known.
It will be used in this form in the algorithm analysis of Section \ref{sec:main-results}. 
Analogous structural assumptions can be found in \eg \citep{Richards2018} to match desired properties of the parametric approximant. 
However, in principle, it may be relaxed to approximate structure matching without essential changes to the forthcoming analyses, and so is omitted for simplicity and brevity.

The actor-critic routine is suggested as follows.
Given $0 < r < R$, consider the following optimization problem $\mathcal{AC}(x_k; R, r)$ at the state $x_k, k \in \N_0$, assuming $x(0) \in \B_R$:
\vspace*{-0.9em}
\begin{subequations} \label{eq:actor-critic-stab}
	\begin{align}
	\hspace*{-6pt}\min_{(u,\theta) \in \U \times \Theta } \quad &\mathbb{J}(x_k,u,\theta)  \label{eq:actor-critic-stab-obj} \tag{A-C-Obj}\\
	\text{\sut} \quad \quad  &\hat{J}(x_k,\theta) \leq \hat{J}(x_k,\theta_{k-1}) + \eps_1  \label{eq:actor-critic-stab-1} \tag{C1}\\
	&V(\hat{x}_{k+1}) \leq \hat{J}(\hat{x}_{k+1}^{u},\theta) + \eps_2 \label{eq:actor-critic-stab-2} \tag{C2}\\
	&\hat{J}(\hat{x}_{k+1}^{u},\theta) - \hat{J}(x_k,\theta) \leq - \frac{\delta}{2}w(x_k) + \eps_3 \label{eq:actor-critic-stab-3} \tag{C3}\\
	&q_1(\|x_k\|) \leq \hat{J}(x_k,\theta) \leq q_2(\|x_k\|). \label{eq:actor-critic-stab-4} \tag{C4}
	\end{align}
\end{subequations}
Here, the cost function is $\mathbb{J}:\R^n \times \U \times \R^p \ra \R$, $\eps_{1,2,3} \geq 0$, $q_{1,2} \in \mathcal{K}_{\infty}, \Theta \subset \R^p$ will be described down below, and $\hat{x}_{k+1}^{u}$ is the state prediction at $t=k \delta$ which is done via
\begin{align} \label{eq:state-pred}
	\hat{x}_{k+1}^{u} = x_k + \delta  f(x_k, u),
\end{align}
although other prediction schemes are possible.
It is due to this state ``prediction'' that $\mathcal{AC}(x_k; R,r)$ 
incorporates relaxation terms $\eps_{1,2,3}$. 
In \eqref{eq:actor-critic-stab-obj}, the objective function $\mathbb{J}$ is set so as to minimize the (squared) Bellman error:
\begin{align} \label{eq:A-C-obj-TD}
\mathbb{J} = \left( \rho(x_k,u) + \hat{J}(\hat{x}_{k+1}^{u},\theta_{k-1}) - \hat{J}(x_k,\theta) \right)^2.
\end{align}
Once the system's state $x_k$ is in $\B_{r^\ast} \subset \B_r$ for some $r^\ast \leq r$, which is referred to as the \textit{core ball}, the setting of $(u_k,\theta_k)$ is arbitrary. 
This is dictated by the fact that the optimization problem may become infeasible in a small vicinity of the origin due to the SH behavior.
Nor is one interested in what happens there as far as SH-setting is concerned.
Various variables in $\mathcal{AC}(x_k; R, r)$ as well as the dependence of the core ball size $r^\ast = r^\ast(R,r)$ on the starting and target balls is now described.
First, regarding the positive-definiteness property of $\hat{J}$ as in \eqref{eq:actor-critic-stab-4}, construct $q_{1,2}$ by specifying bounds on the weight norm as follows.
Let $\underline{\theta},\overline{\theta}$ be \sut $\underline{\theta} \leq \|\theta^{\#}\| \leq \overline{\theta}$ (see Assumption \ref{asm:struct-equiv}) and define 
\begin{align} \label{eq:q_1-and-Theta}
\Theta := \{ \theta \in \R^p: \underline \theta \le \nrm{\theta} \le \overline{\theta} \}, \spc & q_1(\|x\|) := \underline{l}(\|x\|) \cdot \underline{\theta}.
\end{align} 
Next, set
\begin{align*}
\bar{J} := \sup_{\subalign{ x \in \ball_R, \theta \in \Theta }} \; \hat{J}(x,\theta).
\end{align*}
Fix an $\eta_R >0$ and specify $R^\ast > R$ \sut $q_1(R^\ast) \geq \bar{J} + \eta_R$. 
Provided with $R^\ast$, let $L_{\varphi} > 0$ be the local Lipschitz constant of $\varphi$ on $\B_{R^\ast}$ and define 
\begin{align} \label{eq:q_2}
q_2(\|x\|) := \overline{\theta} L_{\varphi} \|x\|. 
\end{align}
Furthermore, let $v^\ast = q_1(r)$ and $r^\ast = q_2^{-1}(\frac{v^\ast}{2})$ (the latter exists since $q_2$ is strictly increasing), which also implies $r^\ast \leq r$. 
At this point, note that for any $\theta \in \Theta$,
\begin{align*}
q_1(\|x\|) \leq \hat{J}(x,\theta) \leq v^\ast \; \Ra \; \|x\| \leq r
\end{align*}
and also 
\begin{align*}
	\frac{v^\ast}{2} \leq \hat{J}(x,\theta) \leq q_2(\|x\|) \; \Ra \; \|x\| \geq r^\ast,
\end{align*}
which relate the value of $\hat{J}$ to the facts that $x \in \B_r$ or $x \not \in \B_{r^\ast(R,r)}$, respectively. 
It can be seen that, among other factors to be detailed later, the bounding functions $q_{1,2}$ contribute to the radius of the target ball.

Finally, call an actor-critic sequence ${(u_k,\theta_k)}_{k \in \N_0}$ \textit{admissible} if, for any $k \in \N_0$, \eqref{eq:actor-critic-stab-1}--\eqref{eq:actor-critic-stab-4} are satisfied along the SH-trajectory of \eqref{eq:sys-SH} as long as $x_k \not \in \ball_{r^\ast(R,r)}$.
A single element of an actor-critic sequence will be called an \emph{actor-critic tuple}.
The following section is devoted to the analysis of the above optimization problem. 

\section{Main Results}\label{sec:main-results}

The following result presents necessary conditions under which $\mathcal{AC}(x_k; R,r), 0 < r < R$ 
yields a practically stabilizing control algorithm. 

\begin{thm} \label{thm:stab}
	Consider the control system \eqref{eq:sys} in the SH-mode \eqref{eq:sys-SH} under the optimization $\mathcal{AC}(x_k; R,r)$. 
	Let $q_{1,2}$ be according to \eqref{eq:q_1-and-Theta} and \eqref{eq:q_2}, and $\bar \delta$ be defined as per Remark \ref{rem:pract-stab-policy} for the radii $0<r^\ast < R^\ast$. 
	Assume that there exists an admissible actor-critic sequence ${(u_k,\theta_k)}_{k \in \N_0}$ for $\mathcal{AC}(x_k; R,r)$ 
	along SH-trajectories of \eqref{eq:sys-SH} with a sampling period $0 < \delta < \bar \delta$,
	and under some $\eps_{1,2,3} \geq 0$.
	Then, there exist $0<\bar \delta_0 \leq \bar \delta, \bar{\eps}_{1,3} >0$ with the following property: if the sampling period $\delta > 0$ satisfies $\delta \le \bar \delta_0$ and $\eps_{1} \leq \bar{\eps}_1, \eps_{3} \leq \bar{\eps}_3$, then the control action sequence $u_k, k \in \N_0$ of the actor-critic sequence ${(u_k,\theta_k)}_{k \in \N_0}$ practically stabilizes the origin of \eqref{eq:sys}.
\end{thm}

\begin{pf}
	First, let $L_f$ be the Lipschitz constant of $f$ on $\B_{R^\ast}$ and define
	\begin{align*}
	\bar{f} := \sup_{\subalign{\hspace{4pt}x &\in \B_{R^\ast} \\ u &\in \U }} \; f(x,u).
	\end{align*}
	Now, let ${(u_k,\theta_k)}_{k \in \N_0}$ be an admissible actor-critic sequence. 
	Consider \eqref{eq:actor-critic-stab-1} at time $k+1$:
	\begin{align*}
	\hat{J}(x_{k+1}^{u_k},\theta_{k+1}) \leq \hat{J}(x_{k+1}^{u_k},\theta_{k}) + \eps_1,
	\end{align*}
	where $x_{k+1}^{u_k}$ is the state after applying $u_k$ at time $t = k \delta$. 
	Using the fact that 
	\begin{align*}
	\norm{x_{k+1}^u - \hat{x}_{k+1}^u} 
	\leq \, &\norm{ \int_{k \delta}^{(k+1) \delta} f(x(\tau),u) \, \text{d}\tau - f(x_k,u) \delta  } \\
	\leq \, &\norm{  \int_{k \delta}^{(k+1) \delta} f(x(\tau),u)  - f(x_k,u) \, \text{d}\tau  } \\
	\leq \, &\int_{k \delta}^{(k+1) \delta} \norm{f(x(\tau),u)  - f(x_k,u)} \, \text{d}\tau \\ 
	\leq \, & \int_{k \delta}^{(k+1) \delta} L_f \bar{f} \delta \, \text{d}\tau = L_f \bar{f} \delta^2
	\end{align*}
	for any $k \in \N_0$, $u \in \U$ and $\delta>0$, it holds that
	\begin{align} \label{subeq:Lipschitz-Jhat-xk+1}
	\begin{split}
	\hat{J}(x_{k+1}^{u_k},\theta_{k}) \leq \; &\hat{J}(\hat{x}_{k+1}^{u_k},\theta_{k}) +  L_{\varphi} \|\theta_k\| \norm{x_{k+1}^u - \hat{x}_{k+1}^u} \\
	\leq \, &\hat{J}(\hat{x}_{k+1}^{u_k},\theta_{k}) + \bar{\theta} L_{\varphi} L_f \bar{f} \delta^2,
	\end{split}
	\end{align}
	for any $\theta_k \in \Theta$. 
	Substituting this into \eqref{eq:actor-critic-stab-1} at $k+1$ gives
	\begin{align*}
	\hat{J}(x_{k+1}^{u_k},\theta_{k+1}) \leq \hat{J}(\hat{x}_{k+1}^{u_k},\theta_{k}) + \overline{\theta} L_{\varphi} L_f \bar{f} \delta^2 + \eps_1
	\end{align*}
	and further subtracting $\hat{J}(x_k,\theta_k)$ thereof yields, using \eqref{eq:actor-critic-stab-3},
	\begin{align} \label{subeq:decay-Jhat}
	\begin{split}
	&\hat{J}(x_{k+1}^{u_k},\theta_{k+1}) - \hat{J}(x_k,\theta_k) \\
	& \hspace{5em}\leq \; -\frac{\delta }{2}w(x_k) + \bar{\theta} L_{\varphi} L_f \bar{f} \delta^2 + \eps_1 + \eps_3,
	\end{split}
	\end{align}
	for any $k \in \N_0$.
	
	In the following it is checked whether any trajectory that starts inside the starting ball $x(0) \in \B_R$ is confined to $\B_{R^\ast}$ and converges to $\B_r$. 
	Let $\hat{J}(x_0,\theta_0) \leq \bar{J}$, for any $x(0) \in \B_R$ and $\theta_0 \in \Theta$. 
	Observe that
	\begin{align} \label{subeq:Overshoot-Jhat-1}
	\hat{J}(x^{u_k}(t),\theta_k) \leq \hat{J}(\chi,\theta_k) +  \overline{\theta} L_{\varphi} L_f \bar{f} \delta^2. 
	\end{align}
	for either $\chi = x_k$ or $\chi = x_{k+1}^{u_k}$. 
	Using \eqref{subeq:Lipschitz-Jhat-xk+1} on the right-hand side above and subtracting $\hat{J}(x_{k},w_k)$ on both sides, it holds that
	\begin{align*}
	&\hat{J}(x^{u_k}(t),\theta_k) - \hat{J}(x_{k},\theta_k) \\
	  &\hspace*{5em}\leq \;\hat{J}(\hat{x}_{k+1}^{u_k},\theta_{k})-\hat{J}(x_{k},\theta_k) + 2 \overline{\theta} L_{\varphi} L_f \bar{f} \delta^2,
	\end{align*}
	and furthermore, since the actor-critic tuple $(u_k,\theta_k)$ satisfies \eqref{eq:actor-critic-stab-3},
	\begin{align} \label{subeq:Overshoot-Jhat-2}
	\hat{J}(x^{u_k}(t),\theta_k) \leq \hat{J}(x_{k},\theta_k)\underbrace{-\frac{\delta }{2}w(x_k) + 2 \overline{\theta} L_{\varphi} L_f \bar{f} \delta^2 + \eps_3}_{=:\Delta_k(\delta)}.
	\end{align}
	Then, under consideration of \eqref{subeq:Overshoot-Jhat-1} at $k=0$, $\delta$ need to satisfy $\overline{\theta} L_{\varphi} L_f \bar{f} \delta^2 \leq \eta_R$ as then
	\begin{align*}
	\hat{J}(x_{0},\theta_0) +  \overline{\theta} L_{\varphi} L_f \bar{f} \delta^2 \leq \bar{J} + \overline{\theta} L_{\varphi} L_f \bar{f} \delta^2 \leq q_1(R^\ast),
	\end{align*}
	from which $\hat{J}(x^{u_0}(t),\theta_0) \leq q_1(R^\ast)$ and thus $\|x(t)\| \leq R^\ast$, $t \in [0, \delta)$, follows. 
	However, $\hat{J}(x^{u_k}(t),\theta_k)$ can be upper bounded more strictly as in \eqref{subeq:Overshoot-Jhat-2}, for any $k \in \N_0$, from which the same conclusion follows if $\delta $ is \sut $\Delta_0(\delta) \leq \eta_R$. 
	Then, for all subsequent time steps $k \in \N_0$, boundedness of the SH-trajectory as $\|x(t)\| \leq R^\ast$, $t \in [k \delta, (k+1) \delta)$, follows if the value of $\hat{J}$ is non-increasing sample-wise, which is shown henceforth. 
	Define now the minimal decay rate as
	\begin{align*}
	\bar{w} = \inf_{r^\ast \leq \|x\| \leq R^\ast} \; \frac{w(x)}{2},
	\end{align*}
	by which $\Delta_k(\delta) = \Delta(\delta) $ can be made independent of the current state. 
	In the following, it is shown that the right-hand side of \eqref{subeq:decay-Jhat} is strictly negative for any $k \in \N_0$ until $\hat{J}(x_k,\theta_k) \leq v^\ast$ \ie that $\hat{J}$ decays to some limit sample-wise. 
	Suppose that $\hat{J}(x_k,\theta_k) \geq \frac{v^\ast}{2}$: 
	Under the minimal decay $\bar{w}$, \eqref{subeq:decay-Jhat} reads as
	\begin{align} \label{subeq:decay-Jhat-min}
	\begin{split}
	&\hat{J}(x_{k+1}^{u_k},\theta_{k+1}) - \hat{J}(x_k,\theta_k) \\
	& \hspace{6em}\leq \; -\delta \bar{w} + \overline{\theta} L_{\varphi} L_f \bar{f} \delta^2 + \eps_1 + \eps_3.
	\end{split}
	\end{align}
	Assume, that the right-hand side indeed is strictly negative \ie $\hat{J}(x_{k+1}^{u_k},\theta_{k+1}) < \hat{J}(x_k,\theta_k)$. 
	Then, at some time step $k \in \N$, the state enters the target ball $\B_r$ and furthermore reaches $\frac{v^\ast}{2} \leq \hat{J}(x_k,\theta_k) \leq \frac{3 v^\ast}{4}$. 
	Note that, by the Lipschitz property of the activation function,
	\begin{align*}
	\|x-y\| \leq \frac{1}{\overline{\theta} L_{\varphi}} \eta_r \; \Ra \; | \hat{J}(x,\theta_k) - \hat{J}(y,\theta_k) | \leq \eta_r
	\end{align*}
	for any $\eta_r >0$ and $\theta_k \in \Theta$. 
	Therefore, for $\norm{x^{u_k}(t) - x_k} \leq \frac{1}{\overline{\theta} L_{\varphi}} \eta_r$, it holds that 
	\begin{align*}
	\hat{J}(x^{u_k}(t),\theta_k) \leq \hat{J}(x_k,\theta_k) + \eta_r \leq \frac{3 v^\ast}{4}+ \eta_r. 
	\end{align*}
	
	
	This means, that $\hat{J}(x^{u_k}(t),\theta_k) \leq v^\ast$ and thus $\|x(t)\| \leq r$, $t \in [k \delta,(k+1)\delta)$, if $\delta$ satisfies $L_f \bar{f} \delta \leq \frac{1}{\overline{\theta} L_{\varphi}} \frac{v^\ast}{4}$. 
	Therefore, choosing
	\vspace*{-0.8em}
	\begin{align} \label{subeq:sample-time-stab}
	\hspace*{-0.3em}\bar \delta_0 \leq  \max_{0 < \delta \leq \bar{\delta}}  \left\{ \delta \; \big| \; \overline{\theta} L_{\varphi} L_f \bar{f} \delta \leq  \frac{v^\ast}{4}, \; \overline{\theta} L_{\varphi}L_f \bar{f} \delta^2 \leq \frac{\bar{w}}{10 }\delta \right\}
	\end{align}
	and setting $\bar{\eps}_1 := \frac{\bar{w}}{2} \delta$ and $\bar{\eps}_3 := \frac{3\bar{w}}{10} \delta$, it follows that $x(t) \in \B_{R^\ast}$ due to \eqref{subeq:Overshoot-Jhat-2} \ie all SH-trajectories are bounded, and for $ x_k \in \B_{R^\ast} \setminus \B_{r^\ast}$,
	\begin{align*}
	\hat{J}(x_{k+1}^{u_k},\theta_{k+1}) - \hat{J}(x_k,\theta_k) \leq - \frac{\bar{w}}{10} \delta. 
	\end{align*}
	for any $0< \delta \leq \bar \delta_0$. 
	In \eqref{subeq:sample-time-stab}, $\bar \delta$ represents a sampling time bound for $V$ to have sample-wise decay on $\B_{R^\ast} \setminus \B_{r^\ast}$ (this will be made use of in Theorem \ref{thm:feas}). 
	
	Thus, it can be concluded that there exist a sampling time bound and relaxation terms of the actor-critic optimization problem \sut \eqref{eq:sys} be practically stabilized in the SH-sense \eqref{eq:sys-SH}. 
	The reaching time for the state to enter the target ball $\B_r$ can be determined in a uniform way from the decay rate and the value of $\hat{J}$ \citep[see \eg][]{Clarke1997-stabilization,Osinenko2018}. 
	\qed
\end{pf}

\begin{rem}
	Observe that the sampling time bounds are lower for higher $\overline{\theta}$, which in turn is user-defined. 
	More specifically, in \eqref{eq:q_1-and-Theta}, $\underline{\theta}$, $\overline{\theta}$ are design factors that influence the relation between the sampling time and the overshoot as well as the core ball, chosen \sut $\theta^{\#} \in \Theta$.
\end{rem}

While the previous result states that the system can be practically stabilized if the sampling time as well as the relaxation terms are chosen sufficiently small, it needs to be shown that indeed for all times $k \in \N_0$, there exists an admissible actor-critic sequence $(u_k,\theta_k)_{k \in \N_0}$ satisfying \eqref{eq:actor-critic-stab-2}-\eqref{eq:actor-critic-stab-4}. 
For that matter, recall that Assumption \ref{asm:CLF} ensures the existence of a CLF while Assumption \ref{asm:struct-equiv} establishes a structural richness $\hat{J}$ to capture $V$. 
Therefore, a nominal stabilizing control policy $\mu(\cdot)$ associated with the CLF $V$, along with $\theta^{\#}$, can guarantee existence of admissible actor-critic sequences for $\mathcal{AC}(x_k;R,r)$, as summarized in the following 

\begin{thm} \label{thm:feas}
	Let Assumption \ref{asm:CLF}-\ref{asm:struct-equiv} hold. 
	Let $\eps_1 \leq \bar{\eps}_1$, $\eps_3 \leq \bar{\eps}_3$ in $\mathcal{AC}(x_k; R,r)$ 
	and $0 < \delta \leq \bar \delta_0$, $\theta \in \Theta$ be bounded as per Theorem \ref{thm:stab}. 
	Given $0<r<R$, with the corresponding $0<r^\ast<R^\ast$, $\bar{w}>0$, and $x(0) \in \B_R$, there exists $\underline{\eps}_1 >0$, $\bar{\eps}_2 >0$, $\underline{\eps}_3 >0$, $0<\bar \delta_1 \leq \bar{\delta}_0$ \sut the following holds: if $\underline{\eps}_1 \leq \eps_1 \leq \overline{\eps}_1$, $0 \leq \eps_2 \leq \bar{\eps}_2 $, $\underline{\eps}_3 \leq \eps_3 \leq \overline{\eps}_3$ and $0< \delta \leq \bar \delta_1$, then for all times $k \in \N_0$ where $x_k \not\in \B_{r^\ast}$, there exists an admissible actor-critic tuple giving rise to an admissible actor-critic sequence $(u_k,\theta_k)_{k \in \N_0}$. 
\end{thm}

\begin{pf}
	Consider a current state $x_k \in \B_{R^\ast} \setminus \B_{r^\ast}$ at some time step $k \in \N_0$. 
	By Assumption \ref{asm:struct-equiv}, it holds that
	\begin{enumerate}
		\item[i)] $\hat{J}(x_{k+1}^{\mu(x_k)},\theta^{\#}) - \hat{J}(x_k,\theta^{\#}) \leq - \delta \frac{w(x_k)}{2}$,
		\item[ii)] $ \alpha_1(\|x_k\|) \leq \hat{J}(x_k,\theta^{\#}) \leq \alpha_1(\|x_k\|)$,
	\end{enumerate}
	where i) is attained for $0<\delta \leq \bar{\delta}_0 \leq \bar{\delta}$ (which was satisfied by Theorem \ref{thm:stab}). 
	Given these properties, it needs to be shown that \eqref{eq:actor-critic-stab-1}-\eqref{eq:actor-critic-stab-4} are well posed to mean that these constraints are feasible for all times where $x_k \not\in \B_{r^\ast}$. 
	From i), it follows that
	\begin{align*}
	\hat{J}(\hat{x}_{k+1}^{\mu(x_k)},\theta^{\#}) - \hat{J}(x_k,\theta^{\#}) \leq &- \delta \frac{w(x_k)}{2} + \|\theta^{\#}\| L_{\varphi} L_f \bar{f} \delta^2 \\
	\leq \, &- \delta \bar{w} + \frac{\bar{w}}{10} \delta.
	\end{align*}
	Hence, $\eps_3$ is lower bounded to be at least $\underline{\eps}_3 := \frac{\bar{w}}{10} \delta$, which results in
	\begin{align*}
	\underline{\eps}_3 := \frac{\bar{w}}{10} \delta \leq \eps_3 \leq \frac{3 \bar{w}}{10} \delta =: \overline{\eps}_3. 
	\end{align*}
	Constraint \eqref{eq:actor-critic-stab-2} is satisfied for any $\eps_2 \geq 0$ due to Assumption \ref{asm:struct-equiv}. 
	Next, note that $\theta_k = \theta^{\#}$ satisfies \eqref{eq:actor-critic-stab-1} only if the sum of the value of $\hat{J}$ under $\theta_{k-1}$ and of $\eps_1$ is not less than the value of $V$. 
	Observe, however, that $\hat{J}$ with $\theta_{k-1}$ is lower bounded due to \eqref{eq:actor-critic-stab-2} as in 
	\begin{align*}
	V(\hat{x}) \leq \hat{J}(\hat{x}_{k}^{u_{k-1}},\theta_{k-1}) + \eps_2,
	\end{align*}
	from which it follows that 
	\begin{align*}
	V(\hat{x}_k^{u_{k-1}})  \leq \hat{J}(x_{k}^{u_{k-1}},\theta_{k-1}) + \eps_2 + \overline{\theta} L_{\varphi} L_f \bar{f} \delta^2 
	\end{align*}
	and thus $V(x_k)$ satisfies
	\begin{align} \label{subeq:V-up-bound-feas}
	\begin{split}
	V(x_k^{u_{k-1}}) \leq \;&\hat{J}(x_{k}^{u_{k-1}},\theta_{k-1}) \\
	&\hspace{2em}+ \eps_2 + \underbrace{\overline{\theta} L_{\varphi} L_f \bar{f} \delta^2}_{\leq \frac{\bar{w}}{10} \delta} +  L_V L_f \bar{f} \delta^2. 
	\end{split}
	\end{align}
	Hence it needs to be shown that $\hat{J}(x_k,\theta^{\#}) = V(x_k)$ is feasible in \eqref{eq:actor-critic-stab-1} for $\eps_{1} \leq \bar{\eps}_1$, that arose from the stability requirements, given the fact that \eqref{subeq:V-up-bound-feas} holds. 
	Conversely, since from \eqref{eq:actor-critic-stab-1}
	\begin{align*}
	\hat{J}(x_k^{u_{k-1}},\theta^{\#}) \leq \hat{J}(x_k^{u_{k-1}},\theta_{k-1}) + \eps_1,
	\end{align*}
	$\eps_1$ is lower bounded by the second line of \eqref{subeq:V-up-bound-feas} as 
	\begin{align*}
	0 \leq \eps_2 + \frac{\bar{w}}{10} \delta +  L_V L_f \bar{f} \delta^2 \leq \eps_1 \leq \overline{\eps}_1.
	\end{align*}
	Setting 
	\begin{align*}
	0 \leq \eps_2 \leq \frac{\bar{w}}{10} \delta =: \overline{\eps}_2
	\end{align*}
	and constraining the sample time \sut $0 < \delta \leq \bar \delta_1$, in which  
	\begin{align*}
	\bar{\delta}_1 \leq \max_{0 \leq \delta \leq \bar \delta_0} \; \left\{ \delta \; \big| \; L_V L_f \bar{f} \delta^2 \leq \frac{\bar{w}}{10} \delta \right\},
	\end{align*}
	results in 
	\begin{align*}
	\underline{\eps}_1 := \frac{3\bar{w}}{10} \delta \leq \eps_1 \leq \frac{\bar{w}}{2} \delta =: \overline{\eps}_1, 
	\end{align*}
	under which $\theta_k = \theta^{\#}$ is admissible for \eqref{eq:actor-critic-stab-1}, for any $0 < \delta \leq \bar \delta_1$. 
	Finally, regarding \eqref{eq:actor-critic-stab-4}, $q_{1,2}$ were chosen according to the specified weight norm bounds $\underline{\theta}, \overline{\theta}$, which allowed $\theta^{\#} \in \Theta$. 
	Hence, feasibility of $\mathcal{AC} (x_k;R,r)$ is shown. \qed
\end{pf}

\begin{rem}
	Note that Assumption \ref{asm:CLF}-\ref{asm:struct-equiv} are only necessary in Theorem \ref{thm:feas} in order to serve feasibility in Theorem \ref{thm:stab}. 
\end{rem}
\begin{rem}
	Due to Assumption \ref{asm:struct-equiv}, the Lipschitz constants $L_V$ and $L_{\varphi}$ are related. 
	They may be comprised to $\bar{L} = \max \{L_V,L_{\varphi}\}$ in Theorem \ref{thm:feas}. 
	However, this may lead to a tighter restriction on the sampling time bounds. 
\end{rem}
The above results can be summarized in the following

\begin{thm}
	Consider the control system \eqref{eq:sys} in the SH-mode \eqref{eq:sys-SH} under the actor-critic optimization $\mathcal{AC}(x_k; R,r)$. 
	Let Assumption \ref{asm:CLF}--\ref{asm:struct-equiv} hold and let $\bar \delta$ be as per Remark \ref{rem:pract-stab-policy} for radii $0<r^\ast<R^\ast$. 
	Then, there exist bounds $\underline{\eps}_{1,2,3},\bar{\eps}_{1,2,3}\geq 0$ and $0 < \bar \delta_1 \leq \bar{ \delta}$ \sut for any $x(0) \in \B_R$, there exists an admissible actor-critic sequence, the action sequence of which practically stabilizes \eqref{eq:sys} as per Definition \ref{def:pract-stab}, if $\underline{\eps}_i \leq \eps_i \leq \overline{\eps}_i$, $i=1,2,3$, and $0 < \delta \leq \bar \delta_1$. 
\end{thm}


\section{Case study and discussion}
\label{sec:case-study}

\begin{figure}[h]
	\centering
	\includegraphics[width=0.85\columnwidth]{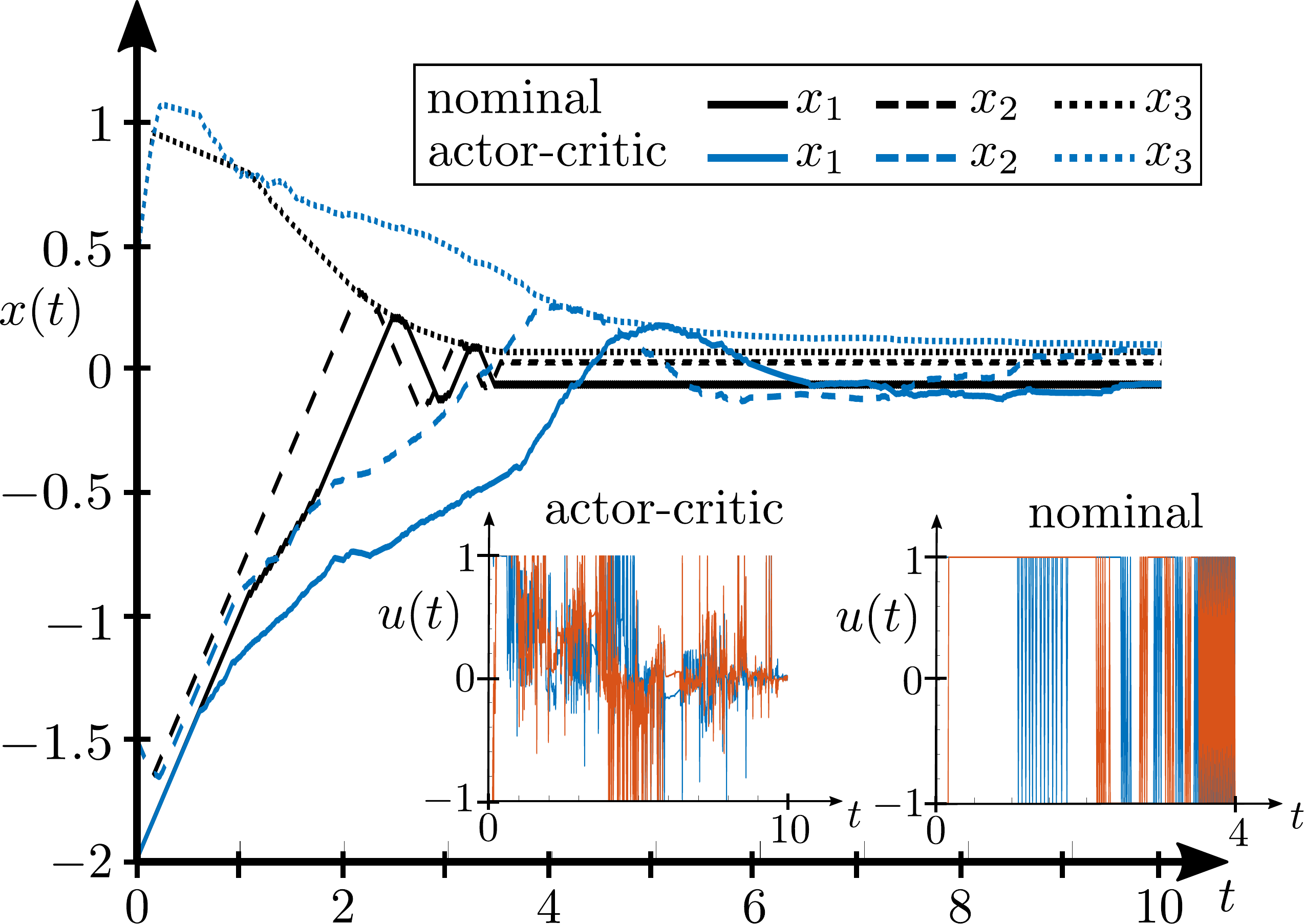}
	\vspace*{-0.5em}
	\caption{State trajectory under the RL-method (blue) and the nominal controller (black) with the corresponding controls ($u_1$ blue, $u_2$ red).
		While a similar, oscillatory-like pattern can be detected in both state trajectories, that under the RL-method converges slower to the target. 
		Yet, the action effort $u$ reaches its constraint boundaries less often compared to the nominal policy. }
	\label{fig:state-input-comb}
\end{figure}

In the following, the non-holonomic integrator
\begin{align*}
\left(\dot{x}_1, \, \dot{x}_2, \,\dot{x}_3 \right)^\top = \left( u_1, \, u_2, \, x_1 u_2 - x_2 u_1 \right)^\top
\end{align*}
is to be practically stabilized under the presented approach. 
It was shown in \eg \citep{Clarke2011}, that the function
\begin{align*}
V(x) = x_1^2 + x_2^2 + 2 x_3^2 - 2 |x_3| \sqrt{x_1^2 +x_2^2 }
\end{align*}
is a global CLF for the system under the input constraint $u \in [-1, \, 1]^2$. 
In the case study, a nominal practically stabilizing control policy is computed via an Inf-convolution technique discussed in \citep{Clarke1997-stabilization,Osinenko2018}.
In the suggested RL-method, the activation function is set to
\begin{align*}
\varphi = \left( x_1^2, \, x_2^2 ,\,2x_3^2, \, - 2 |x_3| \sqrt{x_1^2 +x_2^2 } \right)^\top, 
\end{align*}
so that $\theta^{\#} = \left( 1, \, 1, \, 1, \, 1 \right)^\top$. 
In \eqref{eq:IH-cost} as well as in \eqref{eq:A-C-obj-TD}, $\rho(x,u) = 0.1 x^\top x + 2 u^\top u$. 
First, the constraints \eqref{eq:actor-critic-stab-1}--\eqref{eq:actor-critic-stab-3} are relaxed with $\eps_{1,2,3} = 5\cdot10^{-8}$.
The sampling time is set to $\delta = 0.01 $ and the radius of the target ball is set to $r=0.1$. 
The trajectory of \eqref{eq:sys} under the suggested RL-method in the SH-mode, starting at $x(0) = (-2, \, -1.5, \, 0.4)^\top$ can be seen in Fig.~\ref{fig:state-input-comb}. 

\begin{figure}[h]
	\centering
	\vspace*{-0.5em}
	\includegraphics[width=0.85\columnwidth]{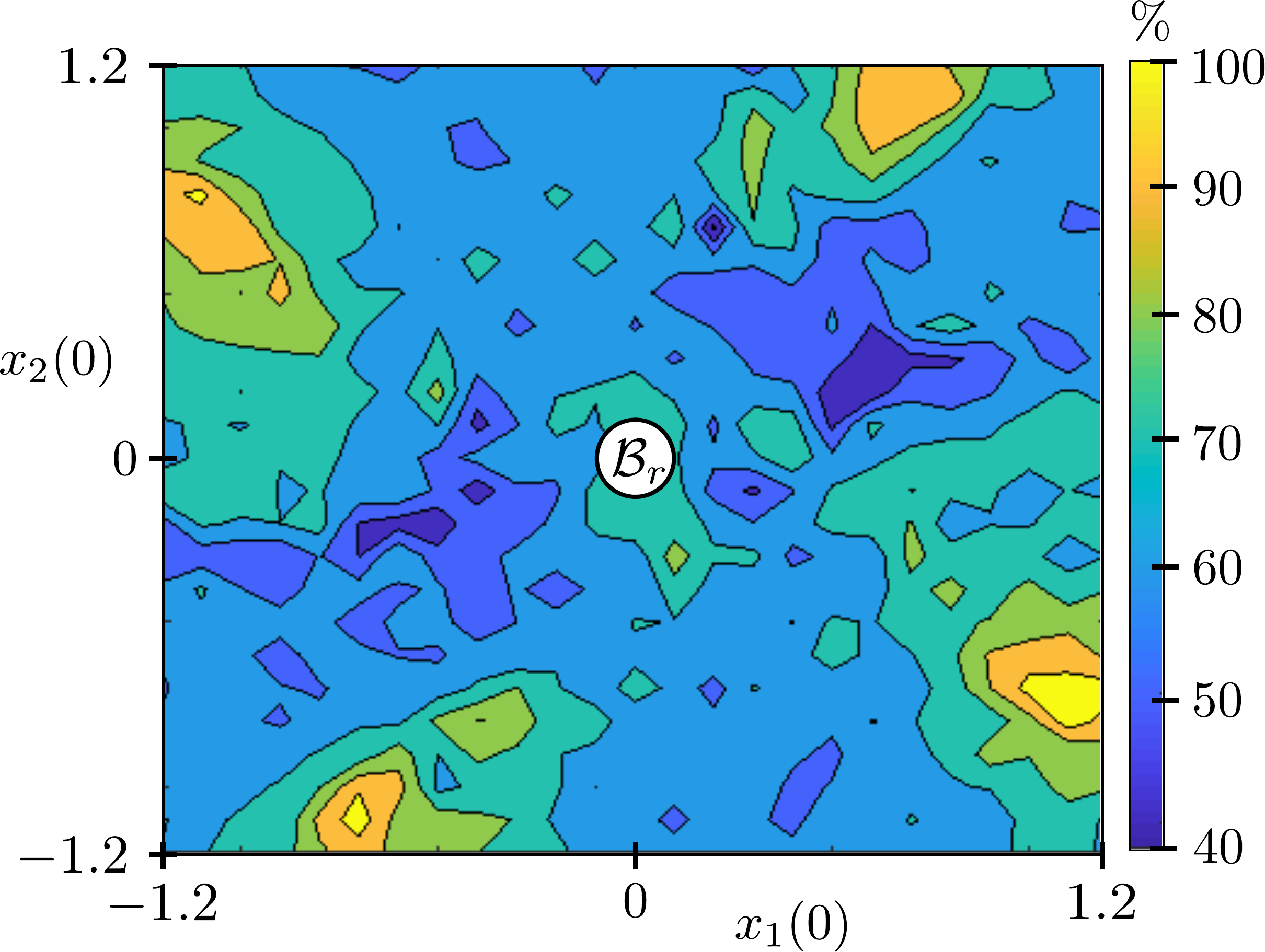}
	\vspace*{-0.5em}
	\caption{Contour of $J_{\text{sim,\%}}(x_0)$. 
	The contour captures the compared cost of driving the state from an initial state in the domain $[-1.2,1.2]^2 \times 0.4$ to a target ball. 
	In most of the grid, the actor-critic control policy could reduce the cost by $20-40\%$. }
	\label{fig:contour-cost}
\end{figure}

From Fig.~\ref{fig:state-input-comb} it can be deduced that stabilization occurs under less action effort whilst allowing the state to converge slower to the target.  
Certain approaches \eg steepest descent, were observed to lead to a bang-bang control chattering between the borders of the set $\U$ \citep[refer to \eg][]{Braun2017,Osinenko2018}, whereas such a behavior could be somewhat alleviated by using the suggested method. 
Consider the simulated, quasi-IH cost
\vspace*{-0.5em}
\begin{align*}
	J_{\text{sim}}[(u_k)_{k \in N_0}](x_0) = \sum_{k=0}^{T-1} \int_{k \delta}^{(k+1)\delta} \rho(x(t),u_k) \, \dt
\end{align*}
under a control sequence $(u_k)_{k \in \N_0}$, where $T$ is the reaching time of the ball $\B_r$ from the starting ball $\B_R$ with $R = 1.75$ \ie the cost of driving the state from $x(0) \in \B_R$ to $\B_r$.   
Using fixed $x_3(0) = 0.4$, the cost difference percentage 
\begin{align*}
	J_{\text{sim,\%}}(x_0) = \frac{J_{\text{sim}}[(u_k)_{k \in N_0}^{\text{actor-critic}}](x_0)}{J_{\text{sim}}[(u_k)_{k \in N_0}^{\text{nominal}}](x_0)} \cdot 100 \%
\end{align*}
on a $(x_1(0),x_2(0)) \in [-1.2, \,1.2]^2 $ grid is depicted in Fig.~\ref{fig:contour-cost}.  

It can be seen that the quasi-IH cost under the RL-method could be improved significantly over a nominal controller. 

\section{Conclusion}
This work was concerned with closed-loop stability issues of RL-methods for dynamical systems.
It suggested to use an initial stabilizability information, specifically, a (non-smooth) CLF, and to introduce it into the constraints of the control scheme.
Practical semi-global stabilizability of the closed-loop, resulting from application of the new method, is analyzed in sample-and-hold manner, which in turn gives insight into the digital realization.
The case study with a non-holonomic integrator showed the merit of the new ideas.

\bibliography{bib/ADP-and-RL,bib/ADP_RL_2,bib/DP,bib/MPC,bib/opt-ctrl,bib/stabilization,bib/MPC_2,bib/nonlin-ctrl,bib/NonlinearControl,bib/intro-bib}

\begin{thebibliography}{23}
\providecommand{\natexlab}[1]{#1}
\providecommand{\url}[1]{\texttt{#1}}
\providecommand{\urlprefix}{URL }
\expandafter\ifx\csname urlstyle\endcsname\relax
  \providecommand{\doi}[1]{doi:\discretionary{}{}{}#1}\else
  \providecommand{\doi}{doi:\discretionary{}{}{}\begingroup
  \urlstyle{rm}\Url}\fi

\bibitem[{Beckenbach et~al.(2018)Beckenbach, Osinenko, and
  Streif}]{Beckenbach2018MPC}
Beckenbach, L., Osinenko, P., and Streif, S. (2018).
\newblock {A}ddressing infinite-horizon optimality in {M}{P}{C} via
  {Q}-learning.
\newblock \emph{IFAC-PapersOnLine}, 51(20), 60--65.

\bibitem[{Bertsekas(2017)}]{Bertsekas2017}
Bertsekas, D.P. (2017).
\newblock Value and policy iteration in optimal control and adaptive dynamic
  programming.
\newblock \emph{IEEE Trans. Neural Netw. Learn. Sys.}, 28(3), 500 -- 509.

\bibitem[{Bian and Jiang(2016)}]{Bian2016}
Bian, T. and Jiang, Z.P. (2016).
\newblock Value iteration and adaptive dynamic programming for data-driven
  adaptive optimal control design.
\newblock \emph{Automatica}, 71, 348--360.

\bibitem[{Braun et~al.(2017)Braun, Gr\"{u}ne, and Kellett}]{Braun2017}
Braun, P., Gr\"{u}ne, L., and Kellett, C.M. (2017).
\newblock Feedback design using nonsmooth control {L}yapunov functions: A
  numerical case study for the nonholonomic integrator.
\newblock In \emph{Proc. of the 56th IEEE Conf. on Decision and Control}.

\bibitem[{Clarke(2011)}]{Clarke2011}
Clarke, F. (2011).
\newblock Lyapunov functions and discontinuous stabilizing feedback.
\newblock \emph{Annual Rev. Control}, 35(1), 13--33.

\bibitem[{Clarke et~al.(1997)Clarke, Ledyaev, Sontag, and
  Subbotin}]{Clarke1997-stabilization}
Clarke, F., Ledyaev, Y., Sontag, E., and Subbotin, A. (1997).
\newblock Asymptotic controllability implies feedback stabilization.
\newblock \emph{IEEE Trans. Automat. Control}, 42(10), 1394--1407.

\bibitem[{Gao and Jiang(2016)}]{Gao2016}
Gao, W. and Jiang, Z.P. (2016).
\newblock {A}daptive {D}ynamic {P}rogramming and adaptive optimal output
  regulation of linear systems.
\newblock \emph{IEEE Trans. Automat. Control}, 61(12), 4164--4169.

\bibitem[{Gao and Jiang(2017)}]{Gao2017}
Gao, W. and Jiang, Z.P. (2017).
\newblock Nonlinear and adaptive suboptimal control of connected vehicles: A
  global {A}daptive {D}ynamic {P}rogramming approach.
\newblock \emph{J. Intelligent \& Robotic Syst.}, 85(3-4), 597--611.

\bibitem[{G{\"o}hrt et~al.(2019)G{\"o}hrt, Osinenko, and
  Streif}]{2019_goeht_adp_lyap}
G{\"o}hrt, T., Osinenko, P., and Streif, S. (2019).
\newblock Adaptive dynamic programming using lyapunov function constraints.
\newblock \emph{IEEE Control Syst. Lett.}, 3(4).

\bibitem[{Heydari(2014)}]{Heydari2014a}
Heydari, A. (2014).
\newblock Revisiting approximate dynamic programming and its convergence.
\newblock \emph{IEEE Trans. on Cyb.}, 44(12), 2733--2743.

\bibitem[{Jiang and Jiang(2015)}]{Jiang2015}
Jiang, Y. and Jiang, Z.P. (2015).
\newblock Global adaptive dynamic programming for continuous-time nonlinear
  systems.
\newblock \emph{IEEE Trans. Automat. Control}, 60(11), 2917--2929.

\bibitem[{{L}iu and {W}ei(2014)}]{Liu2014-PI-ADP-conv-prf}
{L}iu, D. and {W}ei, Q. (2014).
\newblock {P}olicy {I}teration adaptive dynamic programming algorithm for
  discrete-timeime nonlinear systems.
\newblock \emph{IEEE Trans. Neural Netw. Learn. Syst.}, 25(3), 621--634.

\bibitem[{Mu et~al.(2017)Mu, Wang, and He}]{Mu2017}
Mu, C., Wang, D., and He, H. (2017).
\newblock Novel iterative neural dynamic programming for data-based approximate
  optimal control design.
\newblock \emph{Automatica}, 81, 240--252.

\bibitem[{Osinenko et~al.(2018)Osinenko, Beckenbach, and Streif}]{Osinenko2018}
Osinenko, P., Beckenbach, L., and Streif, S. (2018).
\newblock {P}ractical {S}ample-and-{H}old {S}tabilization of {N}onlinear
  {S}ystems {U}nder {A}pproximate {O}ptimizers.
\newblock \emph{IEEE Control Syst. Lett.}, 2(4), 569--574.

\bibitem[{Primbs et~al.(1999)Primbs, Nevisti\'{c}, and Doyle}]{Primbs1999}
Primbs, J.A., Nevisti\'{c}, V., and Doyle, J.C. (1999).
\newblock {N}onlinear optimal control: A control {L}yapunov function and
  receding horizon perspective.
\newblock \emph{Asian J. Control}, 1(1), 14--24.

\bibitem[{Recht(2019)}]{Recht2019}
Recht, B. (2019).
\newblock A tour of reinforcement learning: The view from continuous control.
\newblock \emph{Annual Rev. Control, Robotics and Autonom. Syst.}, 2, 253--279.

\bibitem[{Richards et~al.(2018)Richards, Berkenkamp, and Krause}]{Richards2018}
Richards, S.M., Berkenkamp, F., and Krause, A. (2018).
\newblock The {L}yapunov neural network: {A}daptive stability certification for
  safe learning of dynamical systems.
\newblock Available at arXiv:1808.00924v2 [cs.SY].

\bibitem[{Sokolov et~al.(2015)Sokolov, Kozma, Werbos, and Werbos}]{Sokolov2015}
Sokolov, Y., Kozma, R., Werbos, L.D., and Werbos, P.J. (2015).
\newblock Complete stability analysis of a heuristic approximate dynamic
  programming control design.
\newblock \emph{Automatica}, 59, 9--18.

\bibitem[{Sontag and Sussmann(1995)}]{Sontag1995}
Sontag, E.D. and Sussmann, H.J. (1995).
\newblock Nonsmooth control-lyapunov functions.
\newblock In \emph{Proc. of the 34th IEEE Conf. on Decision and Control}.

\bibitem[{Sutton and Barto(2018)}]{Sutton2018}
Sutton, R.S. and Barto, A.G. (2018).
\newblock \emph{Reinforcement Learning: An Introduction}.
\newblock The MIT Press, 2nd ed. edition.

\bibitem[{Vamvoudakis and Lewis(2010)}]{Vamvoudakis2010}
Vamvoudakis, K.G. and Lewis, F.L. (2010).
\newblock Online actor-critic algorithm to solve the continuous-time infinite
  horizon optimal control problem.
\newblock \emph{Automatica}, 46(5), 878--888.

\bibitem[{Wei et~al.(2016)Wei, Liu, and Lin}]{Wei2016}
Wei, Q., Liu, D., and Lin, Q. (2016).
\newblock {D}iscrete-time local {V}alue {I}teration {A}daptive {D}ynamic
  {P}rogramming: Convergence analysis.
\newblock \emph{IEEE Trans. Syst. Man Cyb.}, 48(6), 875--891.

\bibitem[{Zhang et~al.(2011)Zhang, Cui, and Luo}]{Zhang2011}
Zhang, H., Cui, L., and Luo, Y. (2011).
\newblock {D}ata-driven robust approximate optimal tracking control for unknown
  general nonlinear systems using {A}daptive {D}ynamic {P}rogramming method.
\newblock \emph{IEEE Trans. Neural Netw.}, 22(12), 2226--2236.

\end{thebibliography}

                                                   
\end{document}